\documentclass[11pt,twoside]{article}
\usepackage{amsmath, amsthm, amscd, amsfonts, amssymb, graphicx, color}
\usepackage[bookmarksnumbered, colorlinks]{hyperref} \usepackage{float}
\usepackage{lipsum}
\usepackage{afterpage}
\usepackage[labelfont=bf]{caption}
\usepackage[nottoc,notlof,notlot]{tocbibind} 

\usepackage{lipsum}
\usepackage{fancyhdr}
\newtheorem{theorem}{Theorem}[section]
\newtheorem{lemma}[theorem]{Lemma}

\newtheorem{corollary}[theorem]{Corollary}
\theoremstyle{definition}
\newtheorem{definition}[theorem]{Definition}
\newtheorem{example}[theorem]{Example}

\renewenvironment{proof}{{\bfseries \noindent Proof.}}{~~~~$\square$}
\makeatletter
\def\th@newremark{\th@remark\thm@headfont{\bfseries}}
\makeatletter

\begin{document}
\thispagestyle{plain}
{\noindent Journal of Mathematical Extension \\
Journal Pre-proof}\\
ISSN: 1735-8299\\
URL: http://www.ijmex.com\\
Original Research Paper\\
\vspace*{9mm}

\begin{center}

{\Large \bf 
ON THE RELATIVE $n$-TENSOR NILPOTENT DEGREE OF FINITE GROUPS}

\let\thefootnote\relax\footnote{\scriptsize Received: January 2008; Accepted: February 2009 }

{\bf H. Golmakani}\vspace*{-2mm}\\
\vspace{2mm} {\small  Ferdowsi University of Mashhad} \vspace{2mm}

{\bf  A. Jafarzadeh$^*$\let\thefootnote\relax\footnote{$^*$Corresponding Author}}\vspace*{-2mm}\\
\vspace{2mm} {\small  Quchan University of Technology} \vspace{2mm}

\end{center}

\vspace{4mm}

{\footnotesize
\begin{quotation}
{\noindent \bf Abstract.} In this paper, we generalize the concepts of the relative commutativity degree $d(G,N)$ of a subgroup $N$ of a finite group $G$ and also the tensor degree of a finite group. We introduce the relative $n$-tensor nilpotent degree of a finite group $G$ with respect to a subgroup $H$ of $G$ and investigate some bounds on which.
\end{quotation}
\begin{quotation}

\end{quotation}}

\section{Introduction and Preliminaries}

All groups considered in this paper are finite. Let $G$ be a group with a normal subgroup $N$. Then $(G,N)$ is said to be a pair of groups. Let $G$ and $N$ act on each other and on themselves by conjugation. The nonabelian tensor product $G \otimes N$ is the group generated by the symbols $g \otimes n$ subject to the relations
\begin{eqnarray*}
gg^{'}\otimes n=(^gg^{'}\otimes ^gn)(g\otimes n)
\end{eqnarray*}
\begin{eqnarray*}
g\otimes nn^{'}=(g\otimes n)(^ng\otimes ^nn^{'})
\end{eqnarray*}
for all $g,g^{'}$ in $G$ and $n,n^{'}$ in $N$.
For an element $x\in G$, we consider the tensor centralizer of $x$ as
\begin{eqnarray*}
C_{G}^{\otimes}(x) =\lbrace a \in G :a\otimes x=1_{G\otimes G}\rbrace
\end{eqnarray*}
which is a subgroup of $G$. The intersection of all tensor centralizers of elements of $G$ is called  the tensor center of $G$ and is denoted by $Z^{\otimes}(G)$. The commutator map $\kappa : G \otimes N \longrightarrow [G,N]$ which is given by $g \otimes n \longrightarrow [g,n]$ for all $g\in G$ and $n\in N$, is an epimorphism of groups and we denote $\ker\kappa$ by $J_{2}(G)$. 
We define the tensor upper central series of $G$ as $Z_{1}^{\otimes}(G)=Z^{\otimes}(G)$ and 
\[Z_{n}^{\otimes}(G)=\lbrace {a\in G : [a, x_{1}, \ldots, x_{n-1}]\otimes x_{n}=1 ; \mbox{ for all }x_{1}, \ldots, x_{n}\in G}\rbrace\]
 for all $n\geq2$.
In fact $Z_{n}^{\otimes}(G)/ Z^{\otimes}(G)=Z_{n-1}(G/Z^{\otimes}(G))$
for all $n\geqslant 1$. Hence we have the ascending tensor central series as 
\begin{eqnarray*}
1\leq Z_{1}^{\otimes}(G)=Z^{\otimes}(G)\leq Z_{2}^{\otimes}(G)\leq Z_{3}^{\otimes}(G)\leq \cdots .
\end{eqnarray*}

In \cite{BJR}, the concept of the tensor degree of a group is introduced as
\begin{eqnarray*}
d^{\otimes}(G)= \frac{\vert\lbrace{(x,y)\in G\times G: x\otimes y =1_{\otimes}}\rbrace\vert}{\vert G\vert^{2}}
\end{eqnarray*}
which may be considered as the distance of $G$ from being equal to $Z^{\otimes}(G)$, because $d^{\otimes}(G)=1$ if and only if $Z^{\otimes}(G)=G$. On the other hand, one may easily check that $d^{\otimes}(G)=1$ if and only if $G$ is abelian. 

One of the most important concepts in probabilistic group theory is  the commutativity degree $d(G)$ of a finite group $G$. It is defined in \cite{GR}. Erfanian~et.~ al. in \cite{ERL} generalized the notation of $d(G)$ by defining the relative commutativity degree of a pair of groups $(G,N)$. Let $N$ be a subgroup of $G$. The relative commutativity degree $d(G,N)$ is the probability of commuting an element of $N$ with an element of $G$. It is obviously seen that $d(G)=d(G,G)$ and $d(G,N)=1$ if and only if $N$ is contained in the center of $G$. They also proved the following theorem:
\begin{theorem}(See \cite{ERL} Theorem 3.9)
Let $H$ and $N$ be two subgroups of $G$ such that $N\unlhd G$ and $N\subseteq H$. Then 
\begin{eqnarray*}
d(H,G) \leq d(H/N, G/N)d(N),
\end{eqnarray*}
equality holds if $N\cap [H,G]=1$.
\end{theorem} 

\begin{theorem}\cite{NR2}
Let $G$ be a group and $p$ be the smallest prime divisor of the order of $G$. Then 
\begin{eqnarray*}
\frac{d(G)}{\vert J_{2}(G) \vert} +\frac{\vert Z^{\otimes}(G) \vert}{\vert G \vert}(1-\frac{1}{\vert J_{2}(G) \vert})&\leq  &d^{\otimes}(G)\\&\leq &d(G)-\frac{(p-1)(\vert Z(G) \vert -\vert Z^{\otimes}(G) \vert)}{p\vert G \vert}.
\end{eqnarray*}
\end{theorem}
The special case when $Z^{\otimes}(G)=1$ is described by the next result and has analogies with Theorem 2.8 in \cite{NR1}. There are analogous to the commutativity degree of groups in \cite{BJR,Ellis,NR2,NR3}.
\begin{theorem}\cite{NR2} Let $G$ be a nonabelian group with $Z^{\otimes}(G)=1$ and $p$ be the smallest prime dividing $\vert G \vert$. Then $d^{\otimes}(G) \leq \frac{1}{p}$.
\end{theorem}

\section{Relative n-Tensor Nilpotent Degree of  Groups}

This section is devoted to define the concept of relative $n$-tensor nilpotent degree of a finite group $G$ and a subgroup $H$. Then we obtain some results on this concept. 
\begin{definition}
 Let $H$ be a subgroup of a finite group $G$. We define the relative $n$-tensor nilpotent degree of $H$ in $G$ as 
\begin{eqnarray*}
d_{n}^{\otimes}(H,G)= \frac{\vert\lbrace{(h_{1},\ldots, h_{n},g): [h_{1},\ldots, h_{n}] \otimes g=1_{H\otimes G} , h_{i} \in H, g \in G }\rbrace\vert}{\vert H \vert ^{n} \vert G\vert}.
\end{eqnarray*}
In the special case when $H=G$, it is called the $n$-tensor nilpotent degree of $G$  is denoted by $d_{n}^{\otimes}(G)$.
 \end{definition}
We begin with two following elementary results.
\begin{lemma}\label{lem1}
Let $G$ be a group, $x\in G$ and $H\leq G$, then
\begin{enumerate}
\item[(i)]
 $[H: C^{\otimes}_{G}(x)\cap H] \leq [G: C^{\otimes}_{G}(x)]$;
\item[(ii)]
Equality holds in (i), if $G=HZ^{\otimes}(G)$. The converse is not true.
\end{enumerate}
\end{lemma}
\begin{proof}
(i) Since $C^{\otimes}_{G}(x) \leq G$, we have $HC^{\otimes}_{G}(x)\subseteq G$ and hence
\begin{eqnarray*}
\vert HC^{\otimes}_{G}(x)\vert =\frac{\vert H\vert\vert C^{\otimes}_{G}(x)\vert}{\vert H\cap C^{\otimes}_{G}(x)\vert} \leq \vert G\vert .
\end{eqnarray*}
Therefore
\begin{eqnarray*}
\frac{\vert H\vert}{\vert H\cap C^{\otimes}_{G}(x)\vert} \leq \frac{\vert G\vert}{\vert C^{\otimes}_{G}(x)\vert}.
\end{eqnarray*}
\\
(ii) We know that $Z^{\otimes}(G) =\cap_{x \in G} C^{\otimes}_{G}(x)$. So, if $G = HZ^{\otimes}(G)$, then $G =HC^{\otimes}_{G}(x)$, for all $x\in G$. Thus, 
\begin{eqnarray*}
\vert HC^{\otimes}_{G}(x)\vert =\frac{\vert H\vert\vert C^{\otimes}_{G}(x)\vert}{\vert H\cap C^{\otimes}_{G}(x)\vert} = \vert G\vert .
\end{eqnarray*}
Therefore
\begin{eqnarray}\label{eq1}
[H: C^{\otimes}_{G}(x)\cap H] = [G: C^{\otimes}_{G}(x)],
\end{eqnarray}
 as  required. For the converse, let equation~(\ref{eq1}) holds. Then obviously we have 
 \begin{eqnarray*}
\vert HC^{\otimes}_{G}(x)\vert = \vert G\vert .
\end{eqnarray*}
This does not require to imply  $G=HZ^{\otimes}(G)$. For example, let $G=Q_8=\langle a,b | b^2=a^4=1, b^{-1}ab=a^{-1}\rangle$ and $H=\langle b\rangle$. By Lemma~4.2. of (\cite{NR2}) we have $Z^{\otimes}(G)=1$ and hence $G\neq HZ^{\otimes}(G)$. However, by proof of Theorem~4.3. of (\cite{NR2}) we have  $C^{\otimes}_{G}(a^2)=\langle a \rangle$ and therefore $G=H C^{\otimes}_{G}(a^2)$.
\end{proof}
\begin{theorem}
Let $H\leq G$. Then $d_{n}^{\otimes}(H,G) \leq [G:H]^{n+1} d_{n}^{\otimes}(G)$ for all $n\geq 1$.
\end{theorem}
\begin{proof} Using Lemma~\ref{lem1}, we have
\begin{eqnarray*}   
d_{n}^{\otimes}(H,G) &=& \frac{1}{\vert H\vert ^{n} \vert G \vert} \sum _{x_{1} \in H}\cdots\sum _{x_{n} \in H}\vert C^{\otimes}_{G}([x_{1},\ldots ,x_{n}]) \vert \\ &=& \frac{1}{\vert H\vert ^{n}}\sum _{x_{1} \in H}\cdots\sum _{x_{n} \in H}\frac{\vert C^{\otimes}_{G}([x_{1},\ldots ,x_{n}])\vert}{\vert G\vert} \\&\leq &  \frac{1}{\vert H\vert ^{n} }\sum _{x_{1} \in H}\cdots\sum _{x_{n} \in H}\frac{\vert C^{\otimes}_{G}([x_{1},\ldots ,x_{n}]) \cap H\vert}{\vert H \vert}\\ &\leq & \frac{1}{\vert H\vert ^{n}}\sum _{x_{1} \in H}\cdots\sum _{x_{n} \in H}\frac{\vert C^{\otimes}_{G}([x_{1},\ldots ,x_{n}])\vert}{\vert H \vert} \\ &\leq &\frac{\vert G\vert ^{n+1}}{\vert H\vert ^{n+1} \vert G \vert ^{n+1}}\sum _{x_{1} \in G}\cdots\sum _{x_{n} \in G}\vert C^{\otimes}_{G}([x_{1},\ldots ,x_{n}])\vert \\&=& [G:H]^{n+1} d_{n}^{\otimes}(G)
\end{eqnarray*}
\end{proof}

\begin{theorem}\label{th1}
Let $H \leq G$. Then
\begin{eqnarray*}
d_{n+1}^{\otimes}(H,G) \leq \frac{1}{2}(1+d_{n}^{\otimes}(\frac{H}{H \cap Z^{\otimes}(G)})).
\end{eqnarray*}
\end{theorem}
\begin{proof}
Write $\overline{H}$ for $\frac{H}{H \cap Z^{\otimes}(G)}$ and for each $x\in H$ let $\overline{x}$ stands for $x(H \cap Z^{\otimes}(G))$ as an element of $\overline{H}$. We know that
\begin{align*}
&d_{n+1}^{\otimes}(H,G)=\\
&\frac{1}{\vert H\vert ^{n+1} \vert G \vert} \vert \lbrace{(x_{1},\ldots ,x_{n+1},y)  :[x_{1},\ldots ,x_{n+1}] \otimes y=1_{H\otimes G} , x_{i} \in H, y \in G}\rbrace \vert.
\end{align*}
Therefore
\begin{align*}
&\vert H\vert ^{n+1} \vert G \vert d_{n+1}^{\otimes}(H,G)\\
=&\vert\lbrace {(x_{1},\ldots ,x_{n+1},y) : [x_{1},\ldots ,x_{n+1}] \otimes y=1, x_{i} \in H, y \in G} \rbrace\vert \\
=&\sum_{x_{1} \in H}\cdots\sum _{x_{n+1} \in H}\vert C^{\otimes}_{G}([x_{1},\ldots ,x_{n+1}])\vert \\
=&\sum_{x_{1} \in H}\cdots\sum_{x_{n+1} \in H,\ [x_{1},\ldots ,x_{n+1}]\in H \cap Z^{\otimes}(G)} \vert C^{\otimes}_{G}([x_{1},\ldots ,x_{n+1}])\vert \\
+&\sum_{x_{1} \in H}\cdots\sum _{x_{n+1} \in H,\ [x_{1},\ldots ,x_{n+1}]\notin H \cap Z^{\otimes}(G)} \vert C^{\otimes}_{G}([x_{1},\ldots ,x_{n+1}])\vert.
\end{align*}
On the other hand,
\begin{align*}
&d_{n}^{\otimes}(\overline{H})=\\
&\left(\frac{1}{\vert \overline{H} \vert}\right)^{n+1}\vert\lbrace {(\overline{x_{1}},\ldots ,\overline{x_{n+1}}) : [\overline{x_{1}},\ldots ,\overline{x_{n}}]\otimes \overline{x_{n+1}}=1, \overline{x_{i}}\in \overline{H} }\rbrace\vert= \\
&\left(\frac{\vert H \cap Z^{\otimes}(G)\vert}{\vert H\vert}\right)^{n+1}\frac{\vert\lbrace {(x_{1},\ldots ,x_{n+1}) : [x_{1},\ldots ,x_{n+1}] \in H \cap Z^{\otimes}(G), x_{i} \in H }\rbrace\vert}{\vert H \cap Z^{\otimes}(G) \vert ^{n+1}}
\end{align*}
and we have 
\begin{align*}
&\sum _{x_{1} \in H}\cdots\sum _{x_{n+1} \in H,\ [x_{1},\ldots ,x_{n+1}]\in H \cap Z^{\otimes}(G)} \vert C^{\otimes}_{G}([x_{1},\ldots ,x_{n+1}])\vert =\\
&\vert H \vert ^{n+1} d_{n}^{\otimes}(\overline{H})\vert G \vert . 
\end{align*}
Therefore
\begin{align*}
&\vert H \vert ^{n+1}\vert G \vert d_{n+1}^{\otimes}(H,G) \leq\\ &\vert H \vert ^{n+1} d_{n}^{\otimes}(\overline{H})\vert G \vert 
+(\vert H \vert ^{n+1} - \vert H \vert ^{n+1} d_{n}^{\otimes}(\overline{H}))\frac{\vert G\vert}{2}
\end{align*}
and
\begin{eqnarray*}
\frac{\vert H \vert ^{n+1}\vert G\vert}{2} d_{n}^{\otimes}(\overline{H})+\frac{\vert H \vert ^{n+1}\vert G\vert}{2}&=&\frac{\vert H \vert ^{n+1}\vert G\vert}{2} (1+d_{n}^{\otimes}(\overline{H})). 
\end{eqnarray*}
Hence we have
\begin{eqnarray*}
d_{n+1}^{\otimes}(H,G)\leq \frac{1}{2} (1+ d_{n}^{\otimes}(\overline{H})).  
\end{eqnarray*}
\end{proof}
\section{Tensor Nilpotent Groups}
We are ready to define the concept of tensor nilpotency of a group:
\begin{definition} 
Let $G$ be a group. Then $G$ is called tensor nilpotent if $Z_{n}^{\otimes}(G)=G$ for some $n\geq0$. For a tensor nilpotent group $G$, the smallest $c\geq0$ in which $Z_{c}^{\otimes}(G)=G$ is called the tensor nilpotency class or briefly the tensor class of $G$.
\end{definition}

\begin{theorem}
For a finite group G, we  have 
\begin{eqnarray*}
d_{n+1}^{\otimes}(G)\leq \frac{1}{2^{n}}(2^{n}-1+d^{\otimes}(\frac{G}{Z_{n}^{\otimes}(G)}))
\end{eqnarray*}
for all $n\geq 1$. 
\end{theorem}
\begin{proof}
We khow that 
\begin{eqnarray*}
Z_{n}^{\otimes}(G)/Z^{\otimes}(G)=Z_{n-1}(G/Z^{\otimes}(G))
\end{eqnarray*}
for all $n\geq 1$. We proceed by induction on $n$. For $n=1$, by   using Theorem~\ref{th1}, we have
\begin{eqnarray*}
d_{2}^{\otimes}(G)&\leq &\frac{1}{2} (1+d(\frac{G}{G\cap Z^{\otimes}(G)}))\\ &=& \frac{1}{2}(1+d(\frac{G}{Z^{\otimes}(G)})).
\end{eqnarray*}
Using Theorem~\ref{th1} and the induction we have  
\begin{eqnarray*}
d_{n+1}^{\otimes}(G)&\leq & \frac{1}{2} (1+d_{n}^{\otimes}(\frac{G}{G\cap Z^{\otimes}(G)}))\\&\leq &\frac{1}{2}(1+\frac{1}{2^{n-1}}(2^{n-1}-1+d^{\otimes}(\frac{\frac{G}{Z^{\otimes}(G)}}{Z_{n-1}(\frac{G}{Z^{\otimes}(G)})}))\\&=& \frac{1}{2}(1+\frac{1}{2^{n-1}}(2^{n-1}-1+d^{\otimes}(\frac{G}{Z_{n}^{\otimes}(G)}))\\&=&\frac{1}{2}(\frac{1}{2^{n-1}}(2^{n-1}+2^{n-1}-1+d^{\otimes}(\frac{G}{Z_{n}^{\otimes}(G)})))\\&=&\frac{1}{2^{n}}(2^{n}-1+d^{\otimes}(\frac{G}{Z_{n}^{\otimes}(G)})),
\end{eqnarray*}
as required. 
\end{proof}
\begin{theorem}\label{dntensor}
If $G$ is not a tensor nilpotent group of class at most $n$, then 
\begin{eqnarray*}
d_{n}^{\otimes}(G)\leq \frac{2^{n+2}-3}{2^{n+2}}.
\end{eqnarray*}
\end{theorem}
\begin{proof}
Since $G$ is not a tensor nilpotent group of class at most $n$, $Z_{n}^{\otimes}(G)\neq G$ and $G/ Z_{n-1}^{\otimes}(G)$ is a nonabelian group. We khow that $d^{\otimes}(G)\leq d(G)$, therefore using Theorem 2.2 in \cite{ERL} implies that $d^{\otimes}(\frac{G}{Z_{n-1}^{\otimes}(G)})\leq \frac{5}{8}$. So we have 
\begin{eqnarray*} 
d_{n}^{\otimes}(G)&\leq &\frac{1}{2^{n-1}}(2^{n-1}-1+d^{\otimes}(\frac{G}{Z_{n-1}^{\otimes}(G)}))\\&\leq & \frac{1}{2^{n-1}}(2^{n-1}-1+\frac{5}{8})\\&=&\frac{1}{2^{n-1}}(2^{n-1}-\frac{3}{2^{3}})\\&=&\frac{2^{n+2}-3}{2^{n+2}},
\end{eqnarray*}
as required.
\end{proof}
\begin{lemma}
If $G$ is tensor nilpotent  of class at most $n$, then $G$ is nilpotent of class $n$.  
\end{lemma}
\begin{proof}
We khow that $Z_{n}^{\otimes}(G)\leq Z_{n}(G)$, so the result followes.
\end{proof}
\begin{theorem}
If $G$ is a nontrivial group and $Z(G)=1$, then 
\begin{eqnarray*}
d_{n}^{\otimes}(G)\leq \frac{2^{n}-1}{2^{n}}
\end{eqnarray*}
\end{theorem}
\begin{proof}
We proceed by induction on $n$. Let $n=1$ since $Z(G)=1$, $G$ is not nilpotent and Theorem 3 in \cite{Lescot} implies that $d(G)\leq\frac{1}{2}$. We know that $d^{\otimes}(G)\leq d(G)\leq\frac{1}{2}$. Therefore
\begin{eqnarray*}
d_{n+1}^{\otimes}(G)&\leq &\frac{1}{2^{n}}(2^{n}-1+d^{\otimes}(G))\\&\leq &\frac{1}{2^{n}}(2^{n}-1+\frac{1}{2})\\&=&\frac{1}{2^{n}}(2^{n}-\frac{1}{2})\\&=&\frac{2^{n+1}-1}{2^{n+1}}. 
\end{eqnarray*}
\end{proof}

\begin{theorem}
Let $H$ be a proper subgroup of $G$. Then for all $n\geq 1$, we have 
\begin{enumerate}
\item[(i)] If $H\subseteq Z_{n}^{\otimes}(G)$, then $d_{n}^{\otimes}(H,G)=1$.
\item[(ii)] If $H\nsubseteq Z_{n}^{\otimes}(G)$ and $H/ H\cap Z^{\otimes}(G)$ is tensor nilpotent of class at most $n-1$, then $d_{n}^{\otimes}(H,G)=1$.
\item[(iii)] If $H\nsubseteq Z_{n}^{\otimes}(G)$ and $H/ H\cap Z^{\otimes}(G)$ is not tensor nilpotent of class at most $n-1$, then $d_{n}^{\otimes}(H,G)\leq \frac{2^{n+2}-3}{2^{n+2}}$.
\end{enumerate} 
\end{theorem}
\begin{proof}
(i) If $H\subseteq Z_{n}^{\otimes}(G)$, then $[h_{1}, ..., h_{n}]\otimes x=1$ for all $h_{1},\ldots, h_{n}\in H$ and $x\in G$. So 
\begin{eqnarray*}
d_{n}^{\otimes}(H,G)=\frac{1}{\vert H\vert ^{n}\vert G\vert}\sum_{h_{1}\in H}\cdots\sum_{h_{n}\in H}\vert C^{\otimes}_{G}([h_{1},\ldots,h_{n}])\vert =1.
\end{eqnarray*}
(ii) Since $H/ H\cap Z^{\otimes}(G)$ is a tensor nilpotent group of class at most $n-1$, for all $\overline{h_{1}},\ldots,\overline{h_{n}}$ in $H/ H\cap Z^{\otimes}(G)$ where $\overline{h_{i}}=h_{i}H/ H\cap Z^{\otimes}(G)$, $h_{i}\in H$, $i=1,\ldots, n$ and $[\overline{h_{1}},\ldots,\overline{h_{n-1}}]\otimes \overline{h_{n}}=1$. We know that there exists homomorphism $H/ H\cap Z^{\otimes}(G)\otimes H/ H \cap Z^{\otimes}(G)\longrightarrow (H/ H\cap Z^{\otimes}(G))^{'}$ given $[\overline{h_{1}},\ldots,\overline{h_{n-1}}]\otimes \overline{h_{n}}\longrightarrow [\overline{h_{1}},\ldots,\overline{h_{n}}] $. Since $[\overline{h_{1}},\ldots,\overline{h_{n-1}}]\otimes \overline{h_{n}}=1$, then $[\overline{h_{1}},\ldots,\overline{h_{n}}]=1$. Therefore there exist $h_{1}, \ldots, h_{n}$ in $H$ such that $[h_{1}, \ldots, h_{n}]H\cap Z^{\otimes}(G)= H\cap Z^{\otimes}(G)$, so $[h_{1}, \ldots, h_{n}]\in H\cap Z^{\otimes}(G)$. Therefore for all $x$ in $G$, we have $[h_{1}, \ldots, h_{n}]\otimes x=1$. Hence,  
$C^{\otimes}_{G}([h_{1},\ldots,h_{n}])=G$.
Thus,
\begin{eqnarray*}
\vert H^{n}\vert \vert G\vert d_{n}^{\otimes}(H,G)&=&\vert \{(h_{1},\ldots ,h_{n},x)\in H^{n}\times G :[h_{1},\ldots ,h_{n}]\otimes x=1\}\vert \\&=& \sum_{h_{1}\in H}\cdots\sum_{h_{n}\in H}\vert C^{\otimes}_{G}([h_{1},\ldots ,h_{n}])\vert \\&=& \vert H\vert^{n}\vert G\vert, 
\end{eqnarray*}
and so, $d_{n}^{\otimes}(H,G)=1$.\\
(iii) Since $H/ H\cap Z^{\otimes}(G)$ is not  tensor nilpotent of class at most $n-1$, by Theorem~\ref{dntensor} we have 
\begin{eqnarray*}
d_{n-1}^{\otimes}(\frac{H}{H\cap Z^{\otimes}(G)})\leq \frac{2^{n+1}-3}{2^{n+1}}.
\end{eqnarray*}
Hence
\begin{eqnarray*}
d_{n}^{\otimes}(H,G)&\leq &\frac{1}{2}(1+d_{n}^{\otimes}(\frac{H}{H\cap Z^{\otimes}(G)}))\\ &\leq &\frac{1}{2}(1+\frac{2^{n+1}-3}{2^{n+1}})\\ &\leq &\frac{1}{2}(\frac{2^{n+1}+2^{n+1}-3}{2^{n+1}})\\ &=&\frac{2^{n+2}-3}{2^{n+2}}.  
\end{eqnarray*}
\end{proof}

\begin{theorem}\label{th2}
Let $G$ be finite group, $H$ and $N$ be subgroups of $G$ such that $N\unlhd G$ and $N\subseteq H$. Then 
\begin{eqnarray*}
d_{n}^{\otimes}(H,G)\leq d_{n}^{\otimes}(H/ N, G/ N).
\end{eqnarray*}
\end{theorem}
\begin{proof}
We have
\begin{align*}
&\vert H\vert ^{n}\vert G\vert d_{n}^{\otimes}(H,G)\\
=&\vert \lbrace{(h_{1},\ldots ,h_{n},y) : [h_{1},\ldots ,h_{n}]\otimes y=1, h_{i} \in H, y \in G }\rbrace \vert \\
=&\sum_{h_{1}\in H}\cdots\sum_{h_{n}\in H}\vert C^{\otimes}_{G}([h_{1},\ldots ,h_{n}])\vert \\
=&\sum_{h_{1}\in H}\cdots\sum_{h_{n}\in H}\frac{\vert C^{\otimes}_{G}([h_{1},\ldots ,h_{n}])N\vert \vert C^{\otimes}_{N}([h_{1},\ldots ,h_{n}])\vert}{\vert N\vert}\\
\leq &\sum_{h_{1}\in H}\cdots\sum_{h_{n}\in H}\vert C^{\otimes}_{G/N}([h_{1}N,\ldots ,h_{n}N])\vert \vert C^{\otimes}_{N}([h_{1},\ldots ,h_{n}])\vert \\
=&\sum_{t_{1}\in H/N}\sum_{h_{1}\in H}\cdots\sum_{t_{n}\in H/N}\sum_{h_{n}\in H} \vert C^{\otimes}_{G/N}([t_{1},\ldots ,t_{n}])\vert \vert C^{\otimes}_{N}([h_{1},\ldots,h_{n}])\vert \\
=&\sum_{t_{1}\in H/N}\cdots\sum_{t_{n}\in H/N}\vert C^{\otimes}_{G/N}([t_{1},\ldots , t_{n}])\vert \sum_{h_{1}\in H}\cdots\sum_{h_{n}\in H}\vert C^{\otimes}_{N}([h_{1},\ldots , h_{n}])\vert \\
\leq &\vert N\vert^{n+1} \sum_{t_{1}\in H/N}\cdots\sum_{t_{n}\in H/N}\vert C^{\otimes}_{G/N}([t_{1},\ldots ,t_{n}])\vert \\
=&\vert H/N\vert ^{n}\vert G/N\vert d_{n}^{\otimes}(H/N, G/N)\vert N\vert ^{n+1}\\
=&\vert H\vert ^{n}\vert G\vert d_{n}^{\otimes}(H/N, G/N).
\end{align*} 
Therefore
\begin{eqnarray*}
d_{n}^{\otimes}(H,G) \leq d_{n}^{\otimes}(H/ N, G/ N).
\end{eqnarray*}
\end{proof}
\begin{corollary}
If $N \trianglelefteq G$, then $d_{n}^{\otimes}(G) \leq d_{n}^{\otimes}(G / N)$
\end{corollary}
\begin{proof}
Let $H=G$. Then by using Theorem~\ref{th2}, the result follows . 
\end{proof}
\section{Some Examples}
In this section we compute the relative n-tensor nilpotent group for some groups.
\begin{example}  
For the Symmetric group of degree 3, $S_{3}$,
we have
\begin{eqnarray*}
d_{n}^{\otimes}(S_{3})&\leq &\frac{1}{2}(1+d_{n-1}^{\otimes}(\frac{S_{3}}{Z^{\otimes}(S_{3})}))\\&\leq &\frac{1}{2}(1+d_{n-1}^{\otimes}(S_{3}))\\&\leq &\frac{1}{2}(1+\frac{1}{2^{n-2}}(2^{n-2}-1+d^{\otimes}(\frac{S_{3}}{Z^{\otimes}_{n-2}(S_{3})})))\\&\leq &\frac{2^{n}-1}{2^{n}}.
\end{eqnarray*}
\end{example} 
On the other hand, $Z^{\otimes}(S_{m})=1$. Therefore $d_{n}^{\otimes}(S_{m})\leq \frac{2^{n}-1}{2^{n}}$.
\begin{example} 
Let $G= C_{4}$ and $H= 2C_{4}$.
Then
\begin{eqnarray*}
d_{2}^{\otimes}(2C_{4}, C_{4})&=&\frac{1}{\vert 2C_{4}\vert ^{2}\vert C_{4}\vert}\sum_{h_{1}\in 2C_{4}}\sum_{h_{2}\in 2C_{4}}\vert C^{\otimes}_{C_{4}}([h_{1}, h_{2}])\vert \\ &=&\frac{1}{2^{2}\times 4}\times 16=1.
\end{eqnarray*}
\end{example} 
\begin{example} 
Let $G= D_{8}=<a, b | a^{4}= b^{2}=1, ba=a^{-1}b>$ be the dihedral group of order 8 and $H$ the subgroup generated by $\lbrace{a^{2}, ab}\rbrace$.
Let $n\geq 3$. Since $G$ is a nilpotent group of class 2, we have $\gamma_n(G)=1$. Hence
\begin{eqnarray*}
d_{n}^{\otimes}(H, D_{8})&=&\frac{1}{\vert H\vert ^{n}\vert D_{8}\vert}\sum_{h_{1}\in H}\cdots\sum_{h_{n}\in H}\vert C^{\otimes}_{D_{8}}([h_{1},\ldots ,h_{n}])\vert \\&=&\frac{1}{\vert H\vert ^{n}\times 8}\times \vert H\vert ^{n} =\frac{1}{8}.
\end{eqnarray*}
The same is true if we put $G= Q_{8}=<a, b | a^{4}=1, a^2=b^{2}, ba=a^{-1}b>$, the quaternion group of order 8, and $H=\lbrace{a^{2}, ab}\rbrace$.
\end{example}


\begin{center}

\end{center}

{\small

\noindent{\bf Hanieh Golmakani}

\noindent Department of Pure Mathematics

\noindent Ferdowsi University of Mashhad

\noindent Mashhad, Iran

\noindent E-mail: h.golmakani@mshdiau.ac.ir}\\

{\small
\noindent{\bf  Abbas Jafarzadeh}

\noindent  Department of Mathematics

\noindent Assistant Professor of Mathematics
 
\noindent Quchan University of Technology

\noindent Quchan, Iran

\noindent E-mail: jafarzadeh.a@qiet.ac.ir, abbas.jafarzadeh@gmail.com}\\

\end{document}